\documentclass{crm-eca}

\usepackage{Geomvap1, amscd}  
\usepackage[all,cmtip]{xy}

\newtheorem{theorem}{Theorem}[section]
\newtheorem{lemma}[theorem]{Lemma}

\newtheorem{proposition}[theorem]{Proposition}

\theoremstyle{definition}
\newtheorem{definition}[theorem]{Definition}

\theoremstyle{remark}
\newtheorem{remark}[theorem]{Remark}

\numberwithin{equation}{section}
\newcommand{\ud}{^{(d)}}

\newcommand{\C}{\mathbb C}
\newcommand{\Z}{\mathbb Z}
\newcommand{\Q}{\mathbb Q}
\newcommand{\R}{\mathbb R}

\newcommand{\N}{\mathbb N}

\newcommand{\up}[1]{^{(#1)}}

\newcommand{\wt}{\widetilde}

\newcommand{\cO}{\mathcal O}

\DeclareMathOperator{\Alba}{Alb}

\DeclareMathOperator{\pic}{Pic}

\DeclareMathOperator{\alb}{alb}

\DeclareMathOperator{\vol}{vol}

\begin{document}

\papertitle{The continuous rank function for varieties of maximal Albanese dimension and its applications}


\paperauthor{Lidia Stoppino}
\paperaddress{Universit\`a di Pavia }
\paperemail{lidia.stoppino@unipv.it}



\makepapertitle

\Summary{In this note, I review an aspect of some new techniques introduced recently in collaboration with Miguel \'Angel Barja and Rita Pardini: the construction of the continuous rank function. I give a sketch of how to use this function to prove the  Barja-Clifford-Pardini-Severi inequalities for varieties of maximal Albanese dimension and to obtain the classification of varieties satisfying the equalities.}

\section{Statement of the results}\label{sec: intro and motiv}
We work over $\C$. 
Let $X$ be a smooth projective $n$-dimensional variety and $a\colon X\to A$ a morphism to an abelian $q$-dimensional variety, such that the pullback homomorphism $a^*\colon \pic^0(A)\to \pic^0(X)$ is injective; we call a morphism with such a property {\em strongly generating}. 
The main case to bear in mind is the one when $A=\Alba(X)$ is the  Albanese variety and  $a=\alb_X$ is its Albanese morphism:
in this case $\alb_X^*$ is an isomorphism.
We shall identify $\alpha\in \pic^0(A)$ with $a^*\alpha\in \pic^0(X)$.

Suppose moreover that $X$ is of {\em maximal $a$-dimension}, i.e. that $a$ is finite on its image. 
In particular this implies that $q(X)\geq q=\dim(A)\geq n$, where $q(X)=h^0(X,\Omega^1_X)=\dim (\pic^0(X))$ is the irregularity of $X$.

Let $L$ be any line bundle on $X$. Consider the following integer, which is called {\em the continuous rank of $L$ (\cite[Def. 2.1]{barja-severi})}.
\begin{equation}
h^0_a(X,L):= \min\{h^0(X, L+\alpha),\,\, \alpha \in \pic^0(A)\}.
\end{equation}
\begin{remark}
By semicontinuity, $h^0_a(X, L)$ coincides with $h^0(X, L+\alpha)$ for $\alpha $ general in $\pic^0(A)$, and by Generic Vanishing,  
if $L=K_X+D$ with $D$ nef, then $h^0_a(X, L)=\chi(L)$, the Euler characteristic of $L$. 
\end{remark}
We need also a restricted version of the rank function: for $M\subseteq X$ a smooth subvariety, there exists a non empty open subset of $\pic^0(A)$ such that $h^0(X_{|M}, L+\alpha)$ is constant. We call this value the {\em restricted continuous rank} $h^0_a(X_{|M},L)$.
A first result that highlights the importance of this invariant is the following: 
\begin{proposition}[Barja, \cite{barja-severi}, Thm 3.6.]\label{prop: big}
If $h^0_a(X,L)>0$ then $L$ is big.
\end{proposition}

Recall now that the {\em volume  $\vol_X(L)$} of $L$ (see for instance \cite{lm}) is 
an invariant encoding positivity properties of the line bundle: for example  $\vol_X(L)=L^n$ if $L$ is nef, and $\vol_X(L)>0$ if and only if $L$ is big.

We start with the following general inequalities between the volume of $L$ and its continuous rank.
\begin{theorem}[Barja-Clifford-Pardini-Severi inequalities]\label{thm: BCPS}
The following inequalities hold:
 \begin{itemize}
\item[(i)] $\vol_X (L) \geq n! h^0_a(L)$;
\item[(ii)] If $K_X-L$ is pseudoeffective, then $\vol_X (L) \geq 2n! h^0_a(L)$.
\end{itemize}
\end{theorem}
For the case $n=1$ the inequalities follow from Riemann-Roch and Clifford's Theorem (\cite{articolotto} Lem.~6.13).
For the case $n=2$ and $L=K_X$ inequality (ii) was stated by Severi in 1932, with a wrong proof, and eventually
 proven by Pardini in 2005 \cite{rita-severi}. 
Barja in \cite{barja-severi}, 
proved the inequalities for any $n$ and $L$ nef.
In \cite{articolotto} Barja, Pardini and myself proved the general version Theorem \ref{thm: BCPS} for any line bundle $L$ on $X$, in the form stated above. 
This is done via new techniques introduced in the same paper. Moreover, with our new methods it is possible to solve the problem of {\em classifying} the couples $(X,L)$ that reach the BCPS equalities, obtaining the following general result (\cite{equality}).
\begin{theorem}\label{thm: classification}[\cite{equality}, Thm 1.1, Thm 1.2]
Suppose $h^0_a(X,L)>0$.  
\begin{itemize}
\item[(i)] If $\lambda(L)=n!$, then  $q=n $ and $\deg a=1$ (i.e $a$ is birational).
\item[(ii)] If $K_X-L$ is pseudoeffective  and $\lambda(L)=2n!$, then $q=n$ and  $\deg a=2$. 
\end{itemize}
\end{theorem}
This result was known for $n=2 $ and $L=K_X$ (\cite{BPS}, \cite{LZ2}) but a general classification was out of reach. 

%
In this note I describe in particular one of the techniques of \cite{articolotto}, i.e. the {\em continuous extension of the continuous rank}.
I give an idea of the steps  of the proof of Theorem \ref{thm: BCPS} and of Theorem \ref{thm: classification} that involve the rank function.
Throughout this note, I make assumptions more restrictive than the ones of loc.cit., in order to simplify the exposition. Needless to say, I will hide some technicalities under the carpet.

\section{Continuous rank function}\label{sec: techniques}
\subsection{Set up: Pardini's covering trick}
Let $\mu_d\colon A \to A$ be the multiplication by $d$ on $A$. For  any integer  $d\geq 1$ consider the variety $X\up{d}$ obtained by fibred product as follows:
\begin{equation}\label{diag:-1}
\begin{CD}
\,\,X\up{d}@>\wt{\mu_d}>>X\\
@V{a_d}VV @VV a V\\
A@>{\mu_d}>>A
\end{CD}
\end{equation}
In general, even if we start from $a=\alb_X$, the morphism $a_d$ need not be $\alb_{X\ud}$: 
what is still true is that $a_d$ is strongly generating, as we see from the result below.
\begin{lemma}[\cite{articolotto} Sec. 2.2 and \cite{BPS} Lemma 2.3]\label{lem: covering trick}
The variety $X\ud$ is smooth and connected and the morphism $\wt\mu_d$ is \'etale with the same  monodromy group of $\mu_d$ ($\cong (\Z/d)^{2q}$). We have the following chain of equalities:
$$ 
\ker((a_d\circ\mu_d)^*)=\ker ((a\circ\wt\mu_d)^*)=\pic^0(A)[d]=\ker(\wt\mu_d^*).
$$
In particular, $\ker (a_d^*)=0$.
\end{lemma}


Now, call $L\up{d}:=\wt\mu_d^*(L)$. 
%
Fix $H$ a very ample divisor on $A$; let $M:=a^*H$ and let $M_d$ be a general smooth member of the linear system $|a_d^*H|$. By \cite[Chap.II.8(iv)]{mumford} we have $a_d^*H\equiv d^2H \mod \pic^0(A)$, and hence 
\begin{equation}\label{eq: mum}
M\ud = \wt{\mu_d}^*(a^*H)=a_d^*(\mu_d^*H)\equiv d^2M_d  \mod \pic^0(A).
\end{equation}
\begin{remark}\label{rem: restrizione}
Observe that the assumptions we have on $X$ are verified by $M_d$ for any $d\geq 1$. 
Precisely, the morphism ${a_d}_{|M_d}\colon M_d\to A$ is strongly generating and $M_d$ is of maximal ${a_d}_{|M_d}$-dimension. Moreover, if we have the hypothesis of Theorem \ref{thm: BCPS} (ii), i.e. that $K_X-L$ is pseudoeffective, then $K_{M_d}-L_{|M_d}$ is pseudoeffective.
\end{remark}


\subsection{Continuous continuous rank}
A basic property of the continuous rank with respect to the construction above is the following (see \cite[Prop. 2.8]{barja-severi}):
\begin{equation}\label{eq: mult}
\forall d\in \N\quad \quad  h^0_{a_d}(X\ud , L\ud)=d^{2q}h^0_a(X,L).
\end{equation}
This just follows from the fact that ${\wt{\mu_d}}_*(\cO_{X\ud})=\oplus_{\gamma\in \ker (\mu_d^*)}\gamma$ by Lemma \ref{lem: covering trick}.

\noindent Now we define an extension of the continuous rank for $\R$-divisors of the form $$L_x:=L+xM, \,\,x\in \R.$$ 
We start with the definition over the rationals. 
\begin{definition}
Let  $x\in\mathbb{Q}$, and let $d\in\N$ such that $d^2x=e\in\Z$.
We define
\begin{equation}\label{defi}
h^0_a(X,L_x):=\frac{1}{d^{2q}}h^0_{a_d}(X\ud,L^{(d)}+eM_d).
\end{equation}
\end{definition}
\begin{remark}
Note that by (\ref{eq: mum}) we have that $M_d$ is an integer divisor on $X\ud$ equivalent to $\frac{e}{d^2}M\ud$ modulo $\pic^0(A)$.
For any $k\in \N$, by (\ref{eq: mum}) and (\ref{eq: mult}) we have:
\begin{equation*}
h^0_{a_{dk}}(X^{(dk)}, L^{(dk)}+ ek^2M_{dk})=h^0_{(a_d)_k}((X\ud)^{(k)}, (L^{(d)})^{(k)} + e M_d^{(k)}) =k^{2q}h^0_{a_d}(X\ud, L^{(d)}+eM_d).
\end{equation*}
Using the above equality, it is immediate to see that  given $d,d'\in \N$, $e,e'\in \Z$ such that $\frac{e}{d^2}=x=\frac{e'}{d'^2}$, the formula (\ref{defi}) with $d$ and with $d'$ agrees with the formula with $dd'$, so the definition is independent of the chosen $d$. 
\end{remark}
By studying the properties of this function on $\Q$, we can in particular see that it has the midpoint property, and thus extend it:
\begin{theorem}[\cite{articolotto}, Theorem 4.2]\label{thm: main-rank}
With the above assumptions, the function $h^0_a(X,L_x)$, extends to a continuous  convex  function  $\phi\colon\R\to \R$. For any $x\in \R$ the left derivative has the following form:\begin{equation}\label{eq: derivative}
D^-\phi(x)=\lim_{d\rightarrow \infty}\frac{1}{d^{2q-2}}h^0_{a_d}(X\up{d}_{|M_d},(L_x)^{(d)}), \quad \forall x\in\R.
\end{equation}
\end{theorem}
\begin{remark}\label{rem: volume}
Let us here recall the formula for the derivative of the volume function for $\R$-divisors (see \cite[Cor.C]{lm}).
Fix $x_0:=\max\{ x \mid \vol_X(L_x)=0\}$. There is a continuous extension of the volume function for $\Q$-divisors, $\vol_X(L_x)=\psi(x)\colon \R\to \R$, which is differentiable for $x\ne x_0$  and
\begin{equation}\label{eq: volume}\psi'(x)=\begin{cases}
0 & x<x_0 \\
n\vol_{X|M}(L_x) & x>x_0
\end{cases}\end{equation}
where $\vol_{X|M}(L_x)$ is the {\em restricted volume}.
So, similarly to what happens to the rank function, also the volume extends and the formula for the derivative involves a restricted function. We will soon use this formula.
\end{remark}

\section{Applications}
The power of this new perspective is the following: if we study the BCPS inequalities as a particular case of inequalities between the rank function and the volume function, the proofs become strikingly simple, and we can apply induction via integration.

Now we state the main technical result (see \cite[Sec.2.4]{articolotto}, \cite[sec.2.5]{equality}). 
\begin{lemma}\label{lem: riduzioni}
There exists a $\Q$-divisor $P$ on $X$ such that for any $x\in \R$ and  $d$ high and divisible enough we have:
\begin{align*}
&\vol_{X|M}(L_x) \geq \vol_{X|M}(P_x)=P_x^{n-1}M= \frac{1}{d^{2q}}((P_x)^{(d)})^{n-1}M_d,\\
&\vol_{X\up{d}|{M_d}}({P_x}^{(d)})=((P_x)^{(d)})^{n-1}M_d,\\
&h^0_{a_d}({X\up{d}}_{|{M_d}},(L_x)^{(d)})=h^0_{a_d}({X\up{d}}_{|{M_d}},(P_x)^{(d)})
\end{align*}
\end{lemma}
The key result here is the so-called {\em continuous resolution of the base locus} introduced firstly in \cite[Sec.3]{barja-severi}.


\subsection{BCPS inequalities}\label{sec: BCPS inequalities}
Now we see how the induction step of the proof of the BCPS inequalities ends up in an application of the fundamental theorem of calculus.
We prove here inequality (i) but the proof works exactly in the same way for (ii) (with the right first induction step).
Consider as above the functions $
\psi(x):={\vol}_X(L_x) \,\, \mbox{ and}\  \,\, \phi(x):=h^0_a(X,L_x).$
Using Lemma \ref{lem: riduzioni} and formula (\ref{eq: volume}) we have that 
$$\psi'(x)= \frac{n}{d^{2q}}((P_x)^{(d)})^{n-1}M_d, \quad \quad D^-\phi(x)=\lim_{d\rightarrow \infty}\frac{1}{d^{2q-2}}h^0_{a_d}({X\up{d}}_{|{M_d}},(P_x)^{(d)}).$$
Now, by Remark \ref{rem: restrizione} $M_d$ and ${a_d}_{|M_d}$ satisfy the assumptions, 
and we can prove via the Lemma \ref{lem: riduzioni} that inequality (i) in dimension $n-1$ implies that 
$$
\psi'(x)\geq n! D^-\phi(x)\,\, \mbox{ for any } x\in \R^{\leq 0}.
$$
%
We may thus apply the Fundamental Theorem of Calculus and  compute
\begin{align*}
\vol_X(L)=\psi(0)=\int_{- \infty}^0 \psi'(x)dx&\geq \,n!\,\int_{- \infty}^0D^-\phi(x)dx=\,n!\,\phi(0)=\,n!\,h^0_a(X,L).
\end{align*}

\subsection{Classification of the limit cases}
Both the BCPS inequalities are sharp: we have by Hirzebruch-Riemann-Roch theorem that equality in (i) holds for $X$ an abelian variety and $L$ any nef line bundle on it. 
As for (ii), consider an abelian variety $A$ and a very ample line bundle $N$ on it. Let $B\in |2N|$ a smooth divisor and let $a\colon X\to A$ be the double cover branched along $B$. let $L=a^*(N)$. We have $$\vol_X(L)=2\vol_A(H)=2n!h^0_{id_A}(A, N)=2n!h^0_a(X, L).$$
In Theorem \ref{thm: classification} we see that essentially the cases above are the only ones reaching the equalities. 
Here we give an idea of a step of the proof of (ii). Consider the function
 $$\nu(x):={\vol}_X(L_x)-2n!h^0_a(X,L_x),\,\, x\in \R.$$
One of the key points in the argument in \cite{equality} is to prove that $\nu(x)\equiv 0$ for $x\leq 0$.
We have   $\nu(0)=0$ by assumption.
From Theorem \ref{thm: BCPS} (with some work) we can prove that $\nu(x) \ge 0$ for $x\le 0$.
Hence, it suffices to show that the left derivative $D^-\nu(x)$ is  $\ge 0$  for $x< 0$.
Using Lemma \ref{lem: riduzioni} we have that for any real $x$ smaller or equal than $0$
$$D^-\nu(x)= \lim_{d\rightarrow \infty}\frac{n}{d^{2q-2}}\left(\vol_{M_d}({P_x}^{(d)}_{|M_d})-2(n-1)!h^0_{a_d}({X\up{d}}_{|{M_d}},(P_x)^{(d)})\right).$$
Now we prove that  the right hand expression is greater or equal to 0 using the relative version of Theorem \ref{thm: BCPS} again in dimension $n-1$.


\begin{remark}\label{rem: no clifford}
In Example 7.9 of \cite{articolotto} we proved that for any integer $m\geq 1$ there exist varieties $X_m$ of maximal Albanese dimension such that $\vol(K_{X_m})/\chi(K_{X_m})$ is arbitrarily close to $2n!$ but with Albanese morphism of degree $2^m$, hence far from being a double cover. 
\end{remark}

\begin{remark}
The continuous rank functions can be computed explicitly for abelian varieties, and in some cases for curves (see the Examples of \cite{articolotto}). There are examples where this function is not $\mathcal C^1$ (\cite[Ex.7.3]{articolotto}). 
The regularity properties of these functions, as well as the geometrical meaning of the points of discontinuity of their derivative, still have to be well understood. Some results in this direction can be found in \cite{JP}.
\end{remark}




\end{document}